\documentclass[letterpaper, 10 pt, conference]{ieeeconf}  

\IEEEoverridecommandlockouts                              
\usepackage{color}
\usepackage{amsmath,amssymb}
\newtheorem{theorem}{Theorem}

\newcommand{\feps}{Weakly coupled}
\newcommand{\FEPS}{{weakly coupled}}
\newcommand{\radon}[1]{\mathcal{R}(#1)}

\title{Efficient finite dimensional approximations for the bilinear 
Schr\"odinger equation with bounded variation controls}

\author{
\authorblockN{Nabile Boussa\"{i}d\authorrefmark{1}, Marco 
Caponigro\authorrefmark{2} and Thomas Chambrion\authorrefmark{3}}
\authorblockA{\authorrefmark{1}Laboratoire de math\'ematiques, Universit\'e de 
Franche--Comt\'e, 25030 Besan\c{c}on, France\\
{\tt\small Nabile.Boussaid@univ-fcomte.fr}}
\authorblockA{\authorrefmark{2}\'Equipe M2N,
 Conservatoire National des Arts et M\'etiers, 75003 Paris, France \\
{\tt\small Marco.Caponigro@cnam.fr}}
\authorblockA{\authorrefmark{3}Universit\'e de Lorraine, Institut \'Elie Cartan 
de Lorraine, UMR 7502,  Vandoeuvre-­l\`es-­Nancy,F-­54506, France  \\
CNRS, Institut \'Elie Cartan de Lorraine, UMR 7502, Vandoeuvre-­l\`es-­Nancy, 
F-­54506,
France\\
Inria, CORIDA, Villers-l\`es-Nancy, F-54600, France\\
{\tt\small Thomas.Chambrion@univ-lorraine.fr}}
}
\begin{document}

\maketitle

\begin{abstract}
In the present analysis, we consider the controllability problem of the abstract
Schr\"odinger equation :
\[
\partial_t \psi = A\psi + u B\psi
\]
where $A$ is a skew-adjoint operator, $B$ a control potential and $u$ is the 
control command. 

We are interested by approximation of this equation by finite dimensional 
systems. 

Assuming that $A$ has a pure discrete spectrum and $B$ is in some sense
regular with respect to  $A$ we show that such an approximation is possible. 
More precisely the solutions are approximated by their projections on finite 
dimensional subspaces spanned by the eigenvectors of $A$.

This approximation is uniform in time and in the control $u$, if this control 
has bounded variation with a priori bounded total variation. Hence if these 
finite dimensional systems are controllable with a fixed bound on the total 
variation of $u$ then the system is approximatively controllable. 

The main outcome of our analysis is that we can build solutions for low regular 
controls $u$ such as bounded variation ones and even Radon measures.
\end{abstract}

\section{Extended abstract}

\paragraph{The wellposedness}

Let $\mathcal{H}$ be a separable Hilbert space (possibly infinite dimensional) 
with scalar 
product $\langle \cdot,\cdot\rangle$ and $\|\cdot\|$ the corresponding norm,
$A, B$ be two (possibly unbounded) skew-symmetric operators on $\mathcal{H}$.  
We consider the
formal bilinear control system
\begin{equation}\label{eq:main}
\frac{d}{dt} \psi(t)=A \psi(t) + u(t) B \psi(t),
\end{equation}
where the scalar control $u$ is to be chosen in a set of real functions.

For any real interval $I$, we define
\[
\Delta_{I}:=\{(s,t)\in I^2\mid \, s\leq t\,\}.
\]

{\bf Definition : Propagator on a Hilbert space}\\ 
Let $I$ be a real interval. 
A family $(s,t)\in
\Delta_I\mapsto X(s,t)$ of linear contractions, that is 
Lipschitz maps with Lipschitz constant less than one, on a Hilbert space 
$\mathcal{H}$,
strongly
continuous in $t$ and $s$  and such that
\begin{enumerate}
\item for any $s<r<t$, $X(t,s)=X(t,r)X(r,s)$,
\item $X(t,t)=I_{\mathcal{H}}$,
\end{enumerate}
is called a contraction propagator on $\mathcal{H}$.

Let us now fix some scalar function $u:I\mapsto \mathbf{R}$ and define 
$A(t)=A+u(t)B$. 

Recall that a family $t\in I \mapsto U(t) \in E $, $E$ a subset
of a Banach space $X$, is in $BV(I,E)$
if there exists $N\geq 0$ such that
\[
 \sum_{j=1}^n \|U(t_j)-U(t_{j-1})\|_X \leq N
\]
for any partition $a=t_0<t_1<\ldots<t_n=b$ of the interval $(a,b)$.
The mapping 
\[
U\in BV(I,E) \mapsto \sup_{a=t_0<t_1<\ldots<t_n=b}\sum_{j=1}^n
\|U(t_j)-U(t_{j-1})\|_X
\]
is a semi-norm on $BV(I,E)$ that we
denote with $\|\cdot\|_{BV(I,E)}$.
The semi-norm in $BV(I,E)$ is also called total variation.

{\bf Assumptions }\\
Let $I$ be a real interval and $\mathcal{D}$ dense subset of
$\mathcal{H}$
\begin{enumerate}
\item $A(t)$ is a maximal skew-symmetric operator on $\mathcal{H}$ with domain 
$\mathcal{D}$ ,
\item \label{Ass:BV} $t\mapsto A(t)$ has bounded variation from $I$ to 
$L(\mathcal{D},\mathcal{H} )$, 
where $\mathcal{D}$ is endowed with the graph topology associated to $A(a)$ for 
$a=\inf 
I$\footnote{The bounded variation of $t\mapsto A(t)$ ensures that any choice of 
$s\in I$ will be equivalent.}, 
 \item \label{Ass:ResolventBound}
 $
 M:=\sup_{t\in I} \left\|(1-A(t))^{-1}\right\|_{L(\mathcal{H},\mathcal{D})}< 
\infty,
 $
\end{enumerate}

We do not assume $t \mapsto A(t)$ to be continuous. However as a
consequence of Assumption~\ref{Ass:BV} (see~\cite[Theorem~3]{Edwards})
it admits right and left limit in $L(\mathcal{D},\mathcal{H} )$, $A(t - 
0)=\lim_{\varepsilon
\to 0^{+}}{A(t-\varepsilon)}$, $A(t + 0)=\lim_{\varepsilon
\to 0^{+}}A(t+\varepsilon)$, for all $t \in I$, and $ A(t-0)= A(t+0)$ for all $t
\in I$ except a countable set.

The the core of our analysis is the following result due to Kato 
(see~\cite[Theorem 2 and Theorem 3]{Kato1953}). It gives sufficient conditions 
for the well-posedness of the system~(\ref{eq:main}).
\begin{theorem}\label{thm:kato}
If $t \in I \mapsto A(t)$ satisfies
the above assumptions, then
there exists a unique contraction propagator 
$X:\Delta_{I}\to L(\mathcal{H})$ such that 
if $\psi_0\in \mathcal{D}$ then $X(t,s)\psi_0\in \mathcal{D}$ and for $(t,s)\in 
\Delta_{I}$ 
\[
 \|A(t)X(t,s)\psi_0\|\leq M 
e^{M\|A\|_{BV(I,L(\mathcal{D},\mathcal{H}))}}\|A(s)\psi_0\|.
\]
and in this case $X(t,s)\psi_0$ is strongly left 
differentiable in $t$ and  right differentiable in $s$ with 
derivative (when $t=s$) $A(t+0)\psi_0$ and $-A(t-0)\psi_0$ respectively.

In the case in which $t \mapsto A(t)$ is continuous and
skew-adjoint,
if $\psi_0\in \mathcal{D}$  then $t\in(s,+\infty)\mapsto X(t,s)\psi_0$  is 
strongly continuously differentiable in $\mathcal{H}$ with derivative 
$A(t)X(t,s)\psi_0$.
\end{theorem}
This  theorem addresses the problem of existence of solution for 
the kind of non-autonomous system we consider here under very mild 
assumptions on the control command $u$. For instance we consider bounded 
variation controls. Under some additional assumptions such as the boundedness 
of the control potential, we can also consider Radon measures.

Let us insist on the quantitative aspect of the theorem as it provides an 
estimate on the growth of the solution in the norm associated with $A$. This is 
quantitative aspect is the starting point of the subsequent comments.

\paragraph{Some supplementary regularity}

The regularity of the solution with respect to the natural structure of the 
problem is considered now. In that respect we can adopt two complementary 
strategies : 
\begin{enumerate}
 \item the regularity of the input-output mapping, the flow of the problem, 
is obtained by proving that the control potential if it is regular enough do 
not alter the regularity properties of the uncontrolled problem.
 \item the regularity can be added to the functional setting of the 
wellposedness problem; namely we solve the problem imposing regularity to the 
constructed solution.
\end{enumerate}

For the first strategy we introduce the following definition.

{\bf Definition : \feps }\\
Let $k$ be a non negative real. A couple of \emph{skew-adjoint}
operators $(A,B)$
is \emph{$k$-\FEPS\ } if 
\begin{enumerate}
 \item $A$ is invertible with bounded
inverse from $D(A)$ to $\mathcal{H}$\label{Assumption:FEPSOnA},
 \item for any real $t$, $e^{tB}D(|A|^{k/2})\subset D(|A|^{k/2})$,
 \item there exists $c\geq 0$ and $c'\geq 0$ such that $B-c$ 
and $-B-c'$ generate contraction semigroups on $D(|A|^{k/2})$ for the norm
$\psi 
\mapsto \||A|^{k/2}u\|$.
 \end{enumerate}

We set, for every positive
real $k$,
$$
\|\psi\|_{k/2} = \sqrt{   \langle  | A|^{k} \psi, \psi \rangle}.
$$
The optimal exponential growth is defined by
\[
 c_k(A,B):=\sup_{t\in \mathbf{R}}\frac{\log\|e^{tB}\|_{L(D(|A|^{k/2}),
D(|A|^{k/2})}}{|t|}.
\]

\begin{theorem}\label{PRO_cont_entree_sortie_Hk-BV}
Let $k$ be a non negative real. 
Let $(A,B)$ be $k$-\FEPS\ . 

For any $u\in BV([0,T],\mathbf{R})\cap B_{L^\infty([0,T])}(0,1/\|B\|_A)$, there 
exists 
a family of contraction propagators in $\mathcal{H}$ that extends uniquely as 
contraction propagators
to $D(|A|^{k/2})$ : $\Upsilon^u:\Delta_{[0,T]}\to L(D(|A|^{k/2}))$ such that 
\begin{enumerate}
 \item for any $t\in [0,T]$, for any $\psi_0\in D(|A|^{k/2})$
 \[
  \|\Upsilon_{t}(\psi_0)\|_{k/2}\leq e^{c_k(A,B)\left|\int_0^t u\right|}  
\|\psi_0\|_{k/2}
 \]\label{Contraction}
 \item for any $t\in [0,T]$, for any $\psi_0\in D(|A|^{1+k/2})$ for any $u\in 
BV([0,T],\mathbf{R})\cap B_{L^\infty([0,T])}(0,1/\|B\|_A))$,
 there exists $m$ (depending only on $A$, $B$ and $\|u\|_{L^\infty([0,T])}$)
 \begin{align*}
  \|\Upsilon_{t}(\psi_0)\|_{1+k/2}\leq& 
me^{m\|u\|_{BV([0,T],\mathbf{R})}}\times\\
\times&e^{c_k(A,B)\left|\int_0^t u\right|}  
\|\psi_0\|_{1+k/2}
\end{align*}
\end{enumerate}
Moreover, for every
$\psi_{0}$ in $D(|A|^{k/2})$, the end-point
mapping
\begin{align*}
\Upsilon(\psi_{0}):BV([0,T],K)&\rightarrow  D(|A|^{k/2})\\
u&\mapsto \Upsilon^u(0,T)(\psi_{0})
\end{align*}
is continuous.
\end{theorem}
As announced before these theorem the two strategies were used in complement to 
establish the second point of the theorem.

There is several outcomes to our analysis. First each attainable target from an 
initial state has to be as regular as the initial state and the control 
potential allows it to be. For instance if we consider the harmonic potential 
for the Shr\"odinger equation and we try to control it by a smooth bounded 
potential then from any initial eigenvector one cannot attain non-smooth non 
exponentially decaying states.

\paragraph{A negative result}

An auxiliary result of our analysis is an immediate generalisation of 
the 
famous negative result by Ball, Marsden and Slemrod~\cite{bms} that the 
attainable set is 
included in a countable union of compact sets for the initial Hilbert setting 
of the problem for integrable control laws. We can show that this still holds 
for 
much smaller spaces than the initial Hilbert space, for instance the domains 
and iterated domains of the uncontrolled problem, for much less regular 
controls such as bounded variation function or even Radon measures.
\begin{theorem}\label{Cor:NoExactControllability2}
Let $k$ be a non negative real. Let
$(A,B)$ be $k$-\FEPS\ .  Let $\psi_0\in D(|A|^{k/2})$.Then 
\[
\bigcup_{L,T,a >0}  \left\{ \alpha \Upsilon^{u}(\psi_{0}), 
\|u\|_{BV([0,T],\mathbf{R})} \leq L, t\in [0,T], |\alpha|\leq 
a\right\}
\]
is  a meagre set (in the
sense of Baire) in 
$L^\infty(I,D(|A|^{k/2})$ as a union of relatively compact subsets.
\end{theorem}

\paragraph{Gallerkin approximation for bounded control potential}

In a practical setting it is clear that this negative result is useless. 
Nonetheless we consider that such a negative result has a philosophical 
consequence, natural systems for which regularity of the bounded potential is 
natural cannot be exactly controllable. such an observation is even more 
dramatic if one considers systems with continuous spectrum.

Once such a comment is made, practical questions remains, one of them is the 
question of the approximate controllability (see \cite[Definition 
1]{BVefficiencyCDC}), we choose here to give a 
corollary 
of our analysis that can be helpfull when one want to tackle this issue.


For every Hilbert basis $\Phi=(\phi_k)_{k \in \mathbf{N}}$ of $\mathcal{H}$,
 we define, for every $N$ in $\mathbf{N}$,
$$
\begin{array}{llcl}
 \pi^{\Phi}_N&:\mathcal{H} & \rightarrow & \mathcal{H}\\
      & \psi & \mapsto & \sum_{j\leq N} \langle \phi_j,\psi\rangle \phi_j\,.
\end{array}
$$

{\bf Definition }\\
Let $(A,B)$ be a couple of unbounded operators and
 $\Phi = (\phi_{n})_{n\in \mathbf{N}}$ be an Hilbert basis of $\mathcal{H}$.
Let $N \in \mathbf{N}$ and  denote ${\cal
L}_N^\Phi=\mathrm{span}(\phi_1,\ldots,\phi_N)$.
The \emph{Galerkin approximation} of order $N$ of system \eqref{eq:main}, when 
it makes sense, is the
system
\begin{equation}\label{eq:sigma}
\dot x = (A^{(\Phi,N)} + u B^{(\Phi,N)}) x
\tag{$\Sigma^{\Phi}_{N}$}
\end{equation}
where
$A^{(\Phi,N)}$ and $B^{(\Phi,N)}$, defined by
$$
A^{(\Phi,N)}=\pi^{\Phi}_N A_{\upharpoonright {\mathcal L}_N} \quad \mbox{ and }
\quad
B^{(\Phi,N)}=\pi^{\Phi}_N B_{\upharpoonright{\mathcal L}_N},
$$
are the \emph{compressions} of $A$ and $B$ (respectively) associated with ${\cal
L}_N$.

\bigskip

Notice that if $A$ is skew-adjoint and $\Phi$ an Hilbert basis made of 
eigenvectors of $A$ then
$(A^{(\Phi,N)},B^{(\Phi,N)})$ satisfies the same assumptions as 
$(A,B)$.
We can therefore define the contraction propagator $X^{u}_{(\Phi,N)}(t,0)$ of
\eqref{eq:sigma} associated with bounded variation control $u$. We can also 
write that
$(A^{(\Phi,N)},B^{(\Phi,N)})$ is $k$-\FEPS\ for any prositive real $k$ as the 
weak coupling is actually invariant or at least does not deteriorate by 
compression with respect 
to a basis of eigenvectors of $A$.

In the case of bounded potentials we can state the following proposition. 
Below we consider 
$B$ be bounded in $D(|A|^{k/2})$ which implies that $(A,B)$ is $k$-\FEPS\ .

\begin{theorem}
\label{Prop:GGAforBounded}
Let $k$ be a positive real. Let $A$ with 
domain $D(A)$ be the generator of a contraction semigroup 
and let $B$ be bounded in $\mathcal{H}$ and $D(|A|^{k/2})$ with $B(1-A)^{-1}$ 
compact. Let  $s$ be
non-negative  numbers with
$0\leq s <k$. 
Then
for every $\varepsilon > 0 $, $L\geq 0$, $n\in \mathbf{N}$, and
$(\psi_j)_{1\leq j \leq n}$ in $D(|A|^{k/2})^n$
there exists $N \in \mathbf{N}$
such that
for any $u\in \radon{(0,T]}$,
$$
|u|([0,T])< L \Rightarrow \| \Upsilon^{u}_{t}(\psi_{j}) -
X^{u}_{(N)}(t,0)\pi_{N} \psi_{j}\|_{s/2} < \varepsilon,
$$
for every $t \geq 0$ and $j=1,\ldots,n$.
\end{theorem}

Hence if these 
finite dimensional systems are controllable with a fixed bound on the total 
variation of $u$ then the system is approximatively controllable.

\paragraph{An example Smooth potentials on compact manifolds}

This example motivated the present analysis because of its physical importance. 
We
consider $\Omega$ a compact Riemannian manifold endowed with the associated
Laplace-Beltrami operator $\Delta$ and the associated measure $\mu$, $V,W:\Omega
\to \mathbf{R}$ two smooth functions and the bilinear quantum system 
\begin{equation}\label{EQ-blse-compact-manifold}
 \mathrm{i}\frac{\partial \psi}{\partial t}=\Delta \psi + V \psi +u(t) W \psi.
\end{equation}
With the previous notations, 
$\mathcal{H}=L^2(\Omega,\mathbf{C})$ endowed with the Hilbert product $\langle 
f,g \rangle
=\int_{\Omega} \bar{f} g \mathrm{d}\mu$, $A=-\mathrm{i}(\Delta+V)$ and
$B=-\mathrm{i}W$.
For every $r\geq 0$, $D(|A|^r)=H^r(\Omega,\mathbf{C})$. There exists a Hilbert
basis $(\phi_k)_{k\in \mathbf{N}}$ of $H$ made of eigenvectors of $A$. Each
eigenvalue of $A$ has finite multiplicity. For every $k$, there exists
$\lambda_k$ in $\mathbf{R}$ such that $A\phi_k=\mathrm{i} \lambda_k \phi_k$. The
sequence $(\lambda_k)_k$ tends to $+\infty$ and, up to a reordering, is non
decreasing.  

Since $B$ is bounded from $D(|A|^k)$ to $D(|A|^k)$, $(A,B,\mathbf{R})$ 
satisfies our assumptions and $(A,B)$ is $k$-\FEPS\ for 
non negative real $k$. 

\bibliographystyle{alpha}

\end{document}